\documentclass{amsart}
\usepackage{amssymb,amsmath,amscd}
\usepackage[latin1]{inputenc}
\DeclareMathOperator{\op}{op}

\begin{document}
\newcommand{\hra}{\hookrightarrow}
\newcommand{\lra}{\leftrightarrow}
\newcommand{\la}{\leftarrow}
\newcommand{\ra}{\rightarrow}
\newcommand{\llra}{\longleftrightarrow}
\newcommand{\lla}{\longleftarrow}
\newcommand{\lr}{\longrightarrow}
\newcommand{\Lra}{\Leftrightarrow}
\newcommand{\La}{\Leftarrow}
\newcommand{\Ra}{\Rightarrow}
\newcommand{\Llra}{\Longleftrightarrow}
\newcommand{\Lla}{\Longleftarrow}
\newcommand{\Lr}{\Longrightarrow}

\newcommand{\B}{{\mathbb B}}
\newcommand{\C}{{\mathbb C}}
\newcommand{\E}{{\mathbb E}}
\newcommand{\K}{{\mathbb K}}
\newcommand{\N}{{\mathbb N}}
\newcommand{\Q}{{\mathbb Q}}
\newcommand{\R}{{\mathbb R}}
\newcommand{\T}{{\mathbb T}}
\newcommand{\Z}{{\mathbb Z}}

\newcommand{\cA}{{\mathcal A}}
\newcommand{\cB}{{\mathcal B}}
\newcommand{\cC}{{\mathcal C}}
\newcommand{\cD}{{\mathcal D}}
\newcommand{\cE}{{\mathcal E}}
\newcommand{\cF}{{\mathcal F}}
\newcommand{\cG}{{\mathcal G}}
\newcommand{\cH}{{\mathcal H}}
\newcommand{\cK}{{\mathcal K}}
\newcommand{\cL}{{\mathcal L}}
\newcommand{\cM}{{\mathcal M}}
\newcommand{\cN}{{\mathcal N}}
\newcommand{\cO}{{\mathcal O}}
\newcommand{\cP}{{\mathcal P}}
\newcommand{\cR}{{\mathcal R}}
\newcommand{\cS}{{\mathcal S}}
\newcommand{\cT}{{\mathcal T}}
\newcommand{\cU}{{\mathcal U}}
\newcommand{\cW}{{\mathcal W}}
\newcommand{\cX}{{\mathcal X}}
\newcommand{\cZ}{{\mathcal Z}}

\newcommand{\ga}{\alpha}
\newcommand{\gb}{\beta}
\renewcommand{\gg}{\gamma}
\newcommand{\gG}{\Gamma}
\newcommand{\gd}{\delta}
\newcommand{\eps}{\varepsilon}
\newcommand{\gve}{\varepsilon}
\newcommand{\gk}{\kappa}
\newcommand{\gl}{\lambda}
\newcommand{\gL}{\Lambda}
\newcommand{\go}{\omega}
\newcommand{\gO}{\Omega}
\newcommand{\gvp}{\varphi}
\newcommand{\gvt}{\vartheta}
\newcommand{\gt}{\theta}
\newcommand{\gT}{\Theta}
\newcommand{\gs}{\sigma}

\newcommand{\ol}{\overline}
\newcommand{\ul}{\underline}

\renewcommand{\d}[1]{\partial_{#1}}
\newcommand{\skp}[1]{\langle #1\rangle}
\newcommand{\ds}{\displaystyle}
\newcommand{\eproof}{{~\hfill$ \vartriangleleft$}}
\renewcommand{\Re}{{\rm Re}\,}
\renewcommand{\Im}{{\rm Im}\,}
\newcommand{\forget}[1]{}
\newcommand{\niss}[1]{{#1}} 
\newcommand{\Bew}{\noindent{\em Beweis}}
\newcommand{\Proof}{\noindent{\em Proof}}

 \newtheorem{leer}{\hspace*{-.3em}}[section]

\newenvironment{lemma}[2]%
{\begin{leer}
 \label{#1}
 {\bf Lemma.}      
 {\sl #2}
 \end{leer}}{}

\newenvironment{thm}[2]%
{\begin{leer}
 \label{#1}
 {\bf Theorem.}    
 {\sl #2}
 \end{leer}}{}

\newenvironment{cor}[2]%
{\begin{leer}
 \label{#1}
 {\bf Corollary.}    
 {\sl #2}
 \end{leer}}{}

\newenvironment{dfn}[2]%
{\begin{leer}
 \label{#1}
 {\bf Definition.}     
 {\rm #2 }
 \end{leer}}{}

\newenvironment{ex}[2]%
{\begin{leer}
 \label{#1}
 {\bf Example.}     
 {\rm #2 }
 \end{leer}}{}

\newenvironment{rem}[2]%
{\begin{leer}
 \label{#1}
 {\bf Remark.}     
 {\rm #2 }
 \end{leer}}{}

\newenvironment{extra}[3]%
{\begin{leer}
 \label{#1}
 {\bf #2.}             
 {\rm #3 }
 \end{leer}}{}

\title[$H_\infty$-calculus for Hypoelliptic Operators]{$H_\infty$-calculus for Hypoelliptic Pseudodifferential Operators}
\author{Olesya Bilyj}
\address{Institut für Analysis, Leibniz Universität Hannover, Welfengarten 1, 30167 Hannover,
Germany }
\author{Elmar Schrohe}
\address{Institut für Analysis, Leibniz Universität Hannover, Welfengarten 1, 30167 Hannover,
Germany }
\author{Jörg Seiler}
\address{Department of Mathematical Sciences,
Loughborough University, Leicestershire LE11 3TU, UK}

\begin{abstract}We establish the existence of a bounded
$H_\infty$-calculus for a large class of hypoelliptic
pseudodifferential operators on $\R^n$ and closed manifolds.
\end{abstract}
\maketitle

\section*{Introduction}

Maximal regularity has proven to be a highly efficient concept in the
theory of nonlinear parabolic evolution equations as it can be used
to establish existence and regularity results for nonlinear equations
by studying their linearizations.
 
One way of establishing maximal regularity for a linear evolution equation
$\partial_tu+Au=f$ is to prove that $A$ admits a bounded $H_\infty$-calculus
in the sense of McIntosh \cite{McIn}.
For excellent surveys on maximal regularity and $H_\infty$-calculus see
Denk, Hieber, Pr\"uss
\cite{DHP} or Kunstmann and Weis \cite{KuWe}.

In this short article we will show how a few basic
functional analytic facts about algebras of pseudodifferential operators,
combined with classical techniques developed by Seeley \cite{Seel0} and Kumano-go \cite{K} 
(see also Kumano-go-Tsutsumi  \cite{KT}) 
imply the existence of a bounded $H_\infty$-calculus  for a large class of
sectorially hypoelliptic pseudodifferential operators.

A key point is to focus on the symbols and to establish the calculus on the symbolic
level. This will then imply the existence of a bounded $H_\infty$-calculus for the
associated operators on each Banach space in which zero-order pseudodifferential operators
of the considered type act continuously.

We will first prove our result for symbols in the Hörmander classes $S^m_{\rho,\delta}$
as most readers will be familiar with that calculus.
Then  we will sketch the changes necessary for symbols in the Beals-Fefferman classes
$S^\mu_{\Phi,\gvp}$ and the Weyl-Hörmander classes $S(m,g)$.

Let us recall some basic facts about the $H_\infty$-calculus. 
In the sequel, $\Lambda$ will denote the sector
\begin{equation}\label{h2.B}
 \Lambda=\Lambda(\theta)=\left\{re^{i\varphi}:
 r\ge 0,\;\theta\le\varphi\le 2\pi-\theta\right\},\qquad 0<\theta<\pi,
\end{equation}
in the complex plane.
By $H_\infty$ we denote the space of all bounded holomorphic
functions $f:\C\setminus\Lambda\to\C$, equipped with the supremum norm,
and by $H$ the subspace consisting of  all functions $f$ for which  
$|f(z)|\le c(|z|^d+|z|^{-d})^{-1}$ for suitable $c,d>0$,
depending on $f$.
For every $f\in H_\infty$ there
exists a sequence $(f_j)$ in $H$ converging to $f$ uniformly on
compact subsets of $\C\setminus\Lambda$ such that $\|f_j\|_\infty\le c\|f\|_\infty$
with a constant $c$ independent of $j$ and $f$.
Moreover, each $f\in H_\infty$ has non-tangential boundary values on $\partial \gL$
in $L^\infty$.

Let $A:\cD(A)\subset E\to E$ be a closed and densely defined operator
in a Banach space $E$ such that
\begin{itemize}
\item [(H1)]
$\Lambda\setminus\{0\}$ is contained in the resolvent set of $A$,  
\item [(H2)]$\|\lambda(\lambda-A)^{-1}\|_{\cL(E)}$ is uniformly bounded on
$\Lambda\setminus\{0\}$, and
\item [(H3)]$A$ is injective with dense range.
\end{itemize}
Then  
\begin{equation}\label{h2.A}
  f(A):=\frac{i}{2\pi }\int_{\partial\Lambda} f(\lambda)(A-\lambda)^{-1}\,d\lambda,
  \qquad f\in H,
\end{equation}
defines an element in $\cL(E)$.
Given $f\in H_\infty$, choose a sequence $(f_j)$ in $H$  converging to $f$
as described above.
Then the limit
$$f(A)x=\lim f_j(A)x$$
exists for $x$ in $\cD(A)\cap \text{ran}\,(A)$, which can be shown to be dense in $E$.
The limit is independent of the choice of the sequence and
defines a closable operator $f(A):\cD(A)\to E$.
Its closure is again denoted by $f(A)$.

We say that the
operator $A$ admits a bounded $H_\infty$-calculus with respect to $\C\setminus\Lambda$,
if  $f(A)$ extends to a bounded operator on $E$ and
\begin{equation}\label{h2.C}
 \|f(A)\|_{\cL(E)}\le M\,\|f\|_{\infty} \qquad\text{for all }f\in H_\infty
\end{equation}
with a constant $M$ independent of  $f$.
A particular choice of a bounded holomorphic function is  
$f(z)=z^{it}$, $t\in\R$. It  
implies the boundedness of the purely imaginary powers:
$\|A^{it}\|_{\cL(E)}\le Me^{|t|\theta}$ for all $t\in\R$.
It had been shown before by Dore and Venni \cite{DoVe}, that this
implies maximal regularity, provided $\gt>\pi/2$.
  
In view of the Banach-Steinhaus theorem, it is sufficient to prove estimate \eqref{h2.C} for
$f\in H$. The key to the proof is a thorough understanding of the resolvent of $A$ for estimating 
the operator norm of \eqref{h2.A}.
 
The existence of a bounded $H_\infty$-calculus has been established for many situations.
Amann, Hieber, Simonett \cite{AHS} and Duong, Simonett \cite{DuSi} treat differential
operators on $\R^n$ and on compact manifolds with little regularity in the coefficients;
Duong in \cite{Duon} considers boundary value problems on smooth manifolds,
extending Seeley's work \cite{Seel2} on bounded imaginary powers;
Denk, Dore, Hieber, Pr\"uss, and Venni \cite{DDHPV} investigate
boundary value problems of little regularity;
Abels \cite{Abels} studies the Stokes operator;  
Escher and Seiler \cite{EsSe}  consider
the Dirichlet-Neumann operator for domains of low regularity; Denk, Saal, and Seiler
study Douglis-Nirenberg systems \cite{DSS}. Coriasco, Schrohe, Seiler
show boundedness of the imaginary powers for differential operators
on manifolds with conical singularities \cite{CSS1} and the existence of an $H_\infty$-calculus
for boundary value problems on manifolds with conical singularities \cite{CSS2}.

\section{Functional Analytic Preliminaries}
We recall the concept of $\Psi$-algebras introduced by Gramsch \cite{G}.
\dfn{1.0}{Let $\cA$ be a unital Banach algebra and $\cB$ a Fréchet subalgebra
with a stronger topology and the same unit.
We call $\cB$ a $\Psi$-subalgebra of $\cA$ if
it is spectrally invariant in $\cA$, i.e., if
$\cB^{-1}=\cB\cap \cA^{-1}$. }

\rem{1.5}{A $\Psi$-algebra  has an open group of invertible elements.
According to a result of Waelbroeck \cite{Wa3}, inversion is then continuous.
}

\extra{1.15}{Symbol algebras}{By $S^m_{\rho,\delta}$,
$0\le\delta\le\rho\le1$,
$\delta<1$ we denote the Hörmander class of pseudodifferential
symbols
on $\R^n$ for which the seminorms
$\{q_{\ga,\gb}:\ga,\gb\in\N_0^n\}$  defined by
$$q_{\ga,\gb}(a)=\sup_{x,\xi}|D^\ga_\xi D^\gb_x a(x,\xi)|
\skp\xi^{-m+\rho|\ga|-\gd|\gb|},$$
are finite. This gives $S^m_{\rho,\gd}$ a Fréchet topology.
Moreover, $S^0_{\rho,\delta}$ is a Fréchet algebra with the Leibniz
product $\#$,
associating to two symbols $a$ and $b$ the symbol $a\#b$ of the composition
$\op a\circ\op b$.}

The following theorem is due to R.~Beals \cite{B}, see also Ueberberg \cite{U}.

\thm{1.10}{$S^0_{\rho,\delta}$ is a $\Psi$-subalgebra of $\cL(L^2(\R^n))$. }

\section{The Resolvent as a Pseudodifferential Operator}

\extra{2.0}{General hypoellipticity assumption}{Let
$a\in S^m_{\rho,\delta}$ for some $m\ge0$ and
$0\le\delta<\rho\le1$,
possibly matrix-valued. Assume
that there exist constants $c,C>0$ such that for $x,\xi\in\R^n$, $|\xi|\ge C$,
the spectrum of $a(x,\xi)$ lies outside  $\gL\cup\{|\gl|\le c\}$
and, for $\gl\in\gL$,
\begin{eqnarray}\label{e2.0}
|\partial^\ga_\xi \partial^\gb_x  a(x,\xi)|\ |(a(x,\xi)-\gl)^{-1}|\le
c_{\ga\gb} \skp\xi^{-\rho|\alpha|+\delta|\beta|}.
\end{eqnarray}
}
\rem{2.1}{(a) For $|\xi|\ge C$, $a(x,\xi)$ thus is invertible, and
its spectrum lies in
$$\gO_{x,\xi} =    \{z\in\C\setminus\gL: |z|< 2|a(x,\xi)|\}.$$

(b) Estimate \eqref{e2.0} continues to hold outside $\gO_{x,\xi}$:
For $|\gl|\ge 2|a(x,\xi)|$, we have $|\gl|^{-1}\le (2|a(x,\xi)|)^{-1}\le
|a(x,\xi)^{-1}|/2 $. Hence
\begin{align*}
 |\partial^\alpha_\xi\partial^\beta_xa(x,\xi)||(a(x,\xi)-\lambda)^{-1}|
 &=|\partial^\alpha_\xi\partial^\beta_x a(x,\xi)||\lambda|^{-1}
    \Big|\Big(\frac{a(x,\xi)}{\lambda}-1\Big)^{-1}\Big|\\
 &\le |\partial^\alpha_\xi\partial^\beta_x a(x,\xi)||a(x,\xi)^{-1}|,
 \end{align*}
and the last term can be estimated using \eqref{e2.0} for $\lambda=0$.

(c) Since $(a-\gl)^{-1}=-\gl^{-1}(1-a(a-\gl)^{-1})$ we deduce from  
(b) and the fact that $0$ is not in the spectrum of $a(x,\xi)$, that,
for some $c_0>0$,
\begin{eqnarray}\label{2.1.1}
|(a(x,\xi)-\gl)^{-1}|\le c_0\skp\gl^{-1},\quad |\xi|\ge C, \gl\notin\gO_{x,\xi}.
\end{eqnarray}
}

Following Seeley's classical idea we now construct a parameter-dependent
parame\-trix to $a-\gl$. In fact, a coarse parametrix will be sufficient for
our purposes.

\dfn{2.3}{For $x,\xi\in\R^n,|\xi|\ge C$ and $\gl\notin \gO_{x,\xi}$
we define the sequence $(b_j)$ recursively by $b_0(x,\xi;\gl)=
(a(x,\xi)-\gl)^{-1}$ and
\begin{eqnarray*}
b_{j+1}(x,\xi;\gl)&=&-b_0(x,\xi;\gl)\ \sum_{\substack{|\ga|+k=j+1,\\ 0\le k\le j}}
\frac1{\ga!}
\partial^\ga_\xi  a(x,\xi)D^\ga_x b_k(x,\xi;\gl).
\end{eqnarray*}
Note that $\ga\not=0$ in the summation, so that
$\partial_\xi^\ga(a-\gl)=\partial^\ga_\xi a$.
}

\extra{2.4}{Key observation}{
$\partial^\ga_\xi D^\gb_xb_j$ is a linear combination
of terms of the form
\begin{eqnarray*}\label{e2.4.1}
b_0(x,\xi;\gl)\ \partial^{\ga_1}_\xi D^{\gb_1}_xa(x,\xi)\ b_0 (x,\xi;\gl)\ldots
\partial^{\ga_r}_\xi D_x^{\gb_r}a(x,\xi) b_0(x,\xi;\gl)
\end{eqnarray*}
with suitable $r$ and $\ga_1+\ldots+\ga_r=j+|\ga|$ and
$\gb_1+\ldots+\gb_r=j+|\gb|$.
Indeed this follows from the iteration process together with
the fact that
\begin{eqnarray*}\label{e2.4.2}
\partial_{j} (a(x,\xi)-\gl)^{-1} =
(a(x,\xi)-\gl)^{-1}\ \partial_{j} a(x,\xi)\ (a(x,\xi)-\gl)^{-1}.
\end{eqnarray*}
for an arbitrary derivative $\partial_{j}=\partial_{x_j}$ or $\partial_{j}=\partial_{\xi_j}$.
Note that for $j+|\ga|+|\gb|>0$, we have at least three factors $b_0$ in $\partial^\ga_\xi D^\gb_xb_j$.  }

\dfn{2.8}{We fix a smooth zero-excision function $\gvp$, vanishing for $|\xi|\le C$ and let
for $x,\xi\in\R^n,  \gl\notin\gO_{x,\xi}$,  $N=1,2,\ldots$.
\begin{eqnarray*}
b^N(x,\xi;\gl)&=&\sum_{j<N}\gvp(\xi)b_j(x,\xi;\gl).
\end{eqnarray*}
Moreover, we define the symbols $r^N(\gl)$, $\lambda\in\Lambda$, by
\begin{eqnarray*}
r^N(\gl)&=& (a-\gl)\# b^N(\gl)-1.
\end{eqnarray*}
}

\lemma{2.10}{{\rm (a)} $\gl\mapsto\skp\gl b^N(\gl)$
is bounded continuous from
$\gL$ to $S^{0}_{\rho,\gd}$.
\\
{\rm (b)}
$\gl\mapsto\skp\gl r^N(\gl)$
is a bounded continuous map from
$\gL$ to $S^{m-N(\rho-\gd)}_{\rho,\gd}$.
}

\begin{proof} 
(a) 
Boundedness is an immediate consequence of \ref{2.4} and \eqref{2.1.1};
the resolvent identity
$b_0(\gl)-b_0(\gl_0)=(\gl_0-\gl)b_0(\gl)b_0(\gl_0)$ implies continuity.

For (b) write $r^N(\gl) =
\left((a-\gl)\#b^N(\gl)-q_N(\gl)\right)+(q_N(\gl)-1)$
with  $q_N(\gl)= \sum_{|\ga|<N}\frac1{\ga!}
\partial^\ga_\xi(a-\gl) D^\ga_xb^N(\gl)$. First consider
$$q_N(\gl)-1=\left(\sum_{j+|\ga|<N}+\sum_{j,|\ga|< N, j+|\ga|\ge N}\right)
\partial^\ga_\xi (a-\gl) D^\ga_x(\gvp b_j)(\gl)\ - \ 1.$$
By construction, the first sum equals $\gvp$ and thus differs from  
$1$ by a regularizing symbol.
In the second sum we have $j\not=0$ and $\ga\not=0$. Hence it is a linear combination
of terms with the structure in \ref{2.4}.
Again continuity follows from the resolvent identity, the
boundedness from \eqref{e2.0} and \eqref{2.1.1}.

Finally recall that the difference
$\skp\gl\left((a-\gl)\#b^N(\gl)-q_N(\gl)\right)$ is given by an
oscillatory integral. Its symbol seminorms in $S^{m-N(\rho-\gd)}_{\rho,\gd}$
can be estimated in terms of those for $\partial^\gg_\xi a$ in
$S^{m-N\rho}_{\rho,\gd}$  and those for $\skp\gl D^\gg_xb^N$ in
$S^{N\gd}_{\rho,\gd}$ for $|\gg|=N$.
As both are bounded and the dependence on
$\gl$ is continuous,  the assertion follows.  
\end{proof}

\rem{2.12}{In the same way we can construct
$\tilde b^N(\gl)$ such that
$$\gl\mapsto \skp\gl \tilde r^N(\gl)=
\skp\gl\ \left(\tilde b^N(\gl)\,\#\,(a-\gl)-1\right)$$
is bounded and continuous from $\gL$ to
$S^{m-N(\rho-\gd)}_{\rho,\gd}$.}

\cor{2.14}{Fix $N$ so large that $m-N(\rho-\gd)\le 0$.
Then  
$\ 1+r^N(\gl)$ tends to $1$  in $S^0_{\rho,\gd}$ as $|\gl|\to\infty$.
For large $R$, it is thus invertible with respect
to the Leibniz product on
$\gL_R=\{\gl\in\gL:|\gl|\ge R\}$, since the
group of invertibles is open by \ref{1.5}.
As inversion is continuous, the inverse also tends to $\ 1$;
its seminorms stay bounded.

Repeating the argument with  $\tilde r^N$ we find that
also $a-\gl$ is invertible
on  $\gL_R$ for a possibly larger $R$. Writing $(\cdot)^{-\#}$ for the inverse
with respect to the Leibniz product,
$$(a-\gl)^{-\#}= b^N(\gl)\,\# \, (1+r^N(\gl))^{-\#},\quad\gl\in\gL_R,$$
and its seminorms in $S^0_{\rho,\gd}$ decay like $\skp\gl^{-1}$.

Moreover, the identity
$(1 +r^N(\gl))^{-\#} =1 - r^N(\gl)\,\#\,(1 +r^N(\gl))^{-\#}$
shows that
$$(a-\gl)^{-\#}- b^N(\gl)= s^N(\gl),\quad\gl\in\gL_R,$$
with $\gl\mapsto \skp\gl^2 s^N(\gl)$ bounded and continuous from
$\gL_R$ to $S^{m-N(\rho-\gd)}_{\rho,\gd}$.
}

\cor{2.17}{Replacing $a$ by $a+c$ for some $c>R$ we obtain the
invertibility of $a-\gl$ for all $\gl\in\gL$.
}

\rem{2.15}{In fact, $\gl\to (a-\gl)^{-\#}$ is a holomorphic function on
$\gL_R$ with values in $S^0_{\rho,\gd}$. This follows from the
resolvent identity and the continuity:
\begin{eqnarray*}
\lefteqn{\lim_{\gl\to\gl_0}\frac{(a-\gl)^{-\#}-(a-\gl_0)^{-\#}}{\gl-\gl_0}}\\
& =&
\lim_{\gl\to\gl_0}(a-\gl)^{-\#}\#(a-\gl_0)^{-\#}= (a-\gl_0)^{-\#}\#(a-\gl_0)^{-\#}
.
\end{eqnarray*}  
}

\section{$H_\infty$-calculus}

\extra{3.00}{Functional calculus for symbols}{Let $a$ satisfy the
assumptions in \ref{2.0}. Assume, moreover, that $a-\gl$ is invertible
with respect to the Leibniz product for all $\gl\in\gL$.
As pointed out in Corollary \ref{2.17} this will always be the case after a
suitable shift.
 
For $f\in H$ define the function $f(a)$ on $\R^n\times\R^n$ by
\begin{eqnarray*}
f(a)(x,\xi) = \frac i{2\pi} \int_{\partial\gL} f(\gl) (a-\gl)^{-\#}(x,\xi)\, d\gl.
\end{eqnarray*}
As before,  $(a-\gl)^{-\#}$ denotes the symbol of the Leibniz inverse to $a-\gl$.
The integral converges in $S^0_{\rho,\delta}$ due to the decay property
of $f$ and  since $\gl\mapsto (a-\gl)^{-\#}$
is continuous and decays like $\skp\gl^{-1}$ in all seminorms by Corollary \ref{2.14}.
}

\thm{3.01}{Under these assumptions, $f(a)$ is a symbol in $S^0_{\rho,\delta}$, and
for each symbol seminorm $q$ there is a constant $M_q$, independent of $f$, such that
$q(f(a))\le M_q\ \|f\|_\infty$. }

\begin{proof} 
Let $q$ be the $S^0_{\rho,\gd}$-seminorm  given by
$q(p)=\sup_{x,\xi}|D^\ga_\xi D^\gb_xp(x,\xi)|\skp\xi^{\rho|\ga|-\gd|\gb|}$.
For large $N$ write
$$f(a) =\frac i{2\pi}\int_{\partial\gL} f(\gl)b^N(\gl)\, d\gl +  
\frac i{2\pi}\int_{\partial\gL} f(\gl)s^N(\gl)\, d\gl=:b^N_f+s^N_f.$$
For the first term, we note that $b^N(x,\xi;\gl)$ is a holomorphic function
of $\gl$ away from the spectrum of $a(x,\xi)$ and
$O(\skp\gl^{-1})$ in all seminorms of $S^0_{\rho,\gd}$. 
Moreover, $f(\lambda)$ decays like $O(\skp\gl^{-d})$ for some $d>0$ as $|\gl|\to\infty$. 
For estimating $b_f^N(x,\xi)$ we therefore can replace the contour $\partial \gL$ by $\partial\gO_{x,\xi}$.
Then
\begin{eqnarray*}
\lefteqn{q(b^N_f)\le \frac1{2\pi} \|f\|_\infty}\\
&& \times\sup_{x,\xi}
\left(\text{Length} \
\partial \gO_{x,\xi} \times \sum_{j<N}\sup_{\partial\gO_{x,\xi}}
|D^\ga_\xi D^\gb_x b_j(x,\xi;\gl)|\skp\xi^{\rho|\ga|-\gd|\gb|}\right).
\end{eqnarray*}
The length of the contour is bounded by a constant times
$|a(x,\xi)|$. The specific form of $D^\ga_\xi D^\gb_x b_j$ observed in
\ref{2.4} together with estimate \eqref{e2.0} implies that
\begin{eqnarray*}
|D^\ga_\xi D^\gb_x b_j(x,\xi;\gl)|\skp\xi^{\rho|\ga|-\gd|\gb|}
\le c\ |b_0(x,\xi;\gl)|
\end{eqnarray*}
for a suitable constant $c>0$. Again estimate \eqref{e2.0} shows that
the right hand side is bounded by a constant times $|a(x,\xi)|^{-1}$.
Hence
$q(b^N_f)\le c_q \|f\|_\infty$
with suitable $c_q$.  
 
For the second term, Corollary \ref{2.14} implies that
$$q(s^N_f)\le \frac1{2\pi}\|f\|_\infty \sup_{\gl\in\gL} \skp\gl^2q(s^N(\gl))\
\left|\int_{\partial \gL}\skp\gl^{-2}\, d\gl\right|\le d_q\ \|f\|_\infty$$
for a suitable $d_q$. Together, these two estimates show the assertion.
\end{proof}

\cor{3.02}{Under the above assumptions we have
$$\op (f(a)) =
\frac i{2\pi} \int_{\partial\gL} f(\gl)(\op(a)-\gl)^{-1}\,d\gl,$$
where we consider the pseudodifferential operators on, say, the Schwartz space $\cS$.}

\extra{3.0}{Functional calculus for operators}{Let $E$ be a
Banach space of tempered distributions on $\R^n$ which contains $\cS$
as a dense subspace and for which the mapping
$$\op: S^0_{\rho,\gd}\to \cL(E)$$
is continuous. Let $a$ be a symbol which meets the assumptions
 in \ref{2.0}. Then $\op a:\cS\to\cS$ has a unique closed extension
$A$ in $E$.  We assume that
\begin{itemize}
\item [(H1')] $\gL$ is contained in the resolvent set of $A$.
\item [(H2')] $\|\skp\gl (A-\gl)^{-1}\|_{\mathcal{L}(E)}$ is uniformly bounded in $\lambda\in\Lambda$.
\end{itemize}
We have sharpened condition (H1) (see the introduction) in that we assume that $0$ belongs
to the resolvent set. In this case, (H2') is equivalent to
(H2), and (H3) is automatically fulfilled.  
According to Corollary \ref{2.17},
(H1') and (H2')  will always hold upon replacing $A$
by $A+c$ with a suitably large positive constant $c$.

For $f$ in $H$ we define the operator $f(A)$ by the Dunford integral \eqref{h2.A}.
Assumption (H2') implies that the integral converges in $\cL(E)$
due to  the decay of $f$ near infinity.  }

\thm{3.1}{Under these assumptions
$$\|f(A)\|_{\cL(E)}\le M\ \|f\|_{\infty}$$
for a suitable constant $M$ independent of $f$.
}

\begin{proof} 
As $\gl\to (A-\gl)^{-1} $ is a continuous function
on $\gL$ with values in $\cL(E)$, it suffices to find an
estimate for the integral over $\partial\gL\cap \{|\gl|>R\}$
for large $R$. On this set
\begin{eqnarray*}
(A-\gl)^{-1}=B^N(\gl)+S^N(\gl)
\end{eqnarray*}
with $B^N = \op b^N$ and $S^N=\op s^N$. Then $(A-\gl)^{-1}$, $B^N$ and $S^N$
depend continuously on $\gl$ as operators in $\cL(E)$, and
the norms of $B^N(\gl)$ and $(A-\gl)^{-1}$
in $\cL(E)$ are $O(\skp\gl^{-1})$, that of $S^N(\gl)$ is
$O(\skp\gl^{-2})$. Hence
the norm of
$$\int_{\partial\gL\cap\{|\gl|\ge R\}} f(\gl)S^N(\gl)\, d\gl$$
is bounded by a constant times $\|f\|_\infty$.

In order to estimate the integral involving $B^N(\gl)=\sum_{j<N}\op b_j(\gl)$,
it is sufficient to treat
$$\int_{\partial\gL\cap\{|\gl|\ge R\}}f(\gl)
\op b_j(\gl)\,d\gl=\op\left(\int_{\partial\gL\cap\{|\gl|\ge R\}}f(\gl)
b_j(\gl)\,d\gl\right).$$
In view of the fact that $b_j(x,\xi;\gl)$ is an analytic function of $\gl$ for
$|\gl|>|a(x,\xi)|$ and that $b_j(x,\xi;\gl)$ decays like $\skp\gl^{-1}$ while
$f(\gl)$ decays like  $\skp\gl^{-d}$ for some $d>0$, we can replace the contour
by the contour $\cC_{x,\xi}$ which runs from $Re^{i\gt}$  to $2|a(x,\xi)|e^{i\gt}$
along the straight ray; then clockwise about the origin on a circular arc
to the point $2|a(x,\xi)|e^{i(2\pi-\gt)}$ and then along the ray to
$Re^{i(2\pi-\gt)}$.

Now the same argument as in the proof of Theorem \ref{3.01}
shows that each seminorm for
$\int_{\cC_{x,\xi} }f(\gl)b_j(\gl)\, d\gl$ in the topology
of $S^0_{\rho,\delta}$ is bounded by a multiple of $\|f\|_\infty$.
As this topology is stronger than that of
$\cL(E)$, we obtain that
$$\left\|\int_{\partial\gL\cap\{|\gl|\ge R\}} f(\gl)B^N(\gl)\, d\gl\right\|
\le c\|f\|_\infty.$$
This completes the argument.
\end{proof}

\rem{3.9}{Although many of the arguments can be carried out for general choices
of $\rho$ and $\gd$, the requirement that $\op$ maps $S^0_{\rho,\gd}$ to $\cL(E)$ 
will in general impose rather strict conditions on $\rho$.
If $E$ is an $L^p$-space of Sobolev,
Besov or Triebel-Lizorkin type, for example, we will have
to choose $\rho=1$ whenever $p\not=2$.  
}

\section{The Manifold Case}
{An operator
$A:C^\infty(M,F)\to C^\infty(M,F)$ acting between
sections of a vector bundle $F$ on a smooth closed manifold $M$ is
called a pseudodifferential operator with local symbols in $S^m_{\rho,\gd}$,
if all localizations of $A$ to coordinate neighborhoods are of the form
$\op a$ for a suitable symbol $a\in S^m_{\rho,\delta}$.
One has to assume $1-\rho\le\gd\le \rho$ for the pseudodifferential
calculus on manifolds to make sense.

The pseudodifferential operators with local symbols in $S^0_{\rho,\delta}$ endowed with
the $S^0_{\rho,\gd}$-seminorms on the local symbols then form a Fréchet algebra which is
a $\Psi$-algebra in $\cL(L^2(M,F))$, cf. \cite{Schr1}.

Suppose the local symbols of $A$ satisfy the assumptions of \ref{2.0}.  
We then construct the parameter-dependent parametrix in each coordinate chart.
The operators associated to the local parametrices are
patched to a global parameter-dependent pseudodifferential operator $B^N(\gl)$ on
the manifold.

Similarly to what was shown by Seeley in \cite{Seel0},
the $S^0_{\rho,\delta}$ symbol seminorms for
$B^N(\gl)$ are $O(\skp\gl^{-1})$ and $(A-\gl)B^N(\gl)=1+R^N(\gl)$
for an operator family $R^N$ whose symbol seminorms in $S^0_{\rho,\gd}$ decay like
$\skp\gl^{-1}$ as $\gl\to\infty$. In the same way, a right inverse is obtained,
and the invertibility of $A-\gl$ follows.
Just as above, $1 + R^N(\gl)$ tends to $1$ in the $S^0_{\rho,\gd}$-topology.
It is therefore invertible in $S^0_{\rho,\gd}$ for large $\gl$, and
$$(A-\gl)^{-1}= B^N(\gl)+S^N(\gl)$$
with an operator family $S^N$ whose symbol seminorms in $S^0_{\rho,\gd}$
are $O(\skp\gl^{-2})$ on $\gL$.}

Assuming that $A-\gl$ is invertible for all $\gl$ in $\gL$, we can define
$f(A)$ by the Dunford integral \eqref{h2.A}. The decay of the resolvent
shows that the integral converges in the topology of
$S^0_{\rho,\gd}$.
Moreover, the same analysis as for Theorem \ref{3.01} shows that for each seminorm
$q$ on $S^0_{\rho,\gd}$, we have $q(f(A))\le M_q\|f\|_\infty$.

Next fix  a Banach space $E$ of distributions on $M$ in which $C^\infty(M)$ is dense.
Assume that for each
pseudodifferential operator $P$ on $M$ with local symbols in $S^0_{\rho,\gd}$, the operator
norm of $P$ in $\cL(E)$ can be estimated in terms of the seminorms of its local
symbols.
Then the analysis for Theorem \ref{3.1} shows that
$A$ has a bounded $H_\infty$-calculus.

\section{Symbols of Beals-Fefferman Type}
Let $\Phi,\gvp$ be a pair of weight functions on $\R^n\times \R^n$
in the sense of Beals and Fefferman,
satisfying the standard conditions (equations (1.1) through (1.6) in \cite{Beals}).
Let $\mu\in O(\Phi,\gvp)$, i.e., for suitable constants $c,C$,
\begin{eqnarray*}
&&|\mu(x,\xi)-\mu(y,\eta)|\le C\ \text{if}\ |x-y|\le c\gvp(x,\xi)\ \text{and}\
|\xi-\eta|\le c\Phi(x,\xi)
\\
&&\text{for some real }k,K,m:\quad
c(\gvp\Phi)^{-m}\le e^\mu\gvp^{-k}\Phi^{-K}\le C(\gvp\Phi)^m.
\end{eqnarray*}
Simple examples are the functions $(K,k):=K\log\Phi+k\log\gvp$ for arbitrary  
$K,k\in\R$.
We denote by $S^\mu_{\Phi,\gvp}$ the associated symbol class on $\R^n$, consisting
of all $a=a(x,\xi)$ such that for all multi-indices  $\ga,\gb$,
\begin{eqnarray}
|D^\ga_\xi D^\gb_x a(x,\xi)|\le C_{\ga,\gb}e^\mu\Phi^{-|\ga|}\gvp^{-|\gb|}.
\end{eqnarray}
The associated seminorms give $S^\mu_{\Phi,\gvp}$ a Fréchet topology.
We recover the Hörmander class $S^m_{\rho,\gd}$ for
$\Phi=\skp \xi^\rho$,  $\gvp=\skp\xi^{-\gd}$, and $\mu=\ln\skp\xi^m$.

Each symbol $a\in S^0_{\Phi,\gvp}$ induces a bounded pseudodifferential
operator on $L^2(\R^n)$, and the symbol topology is stronger than the operator
topology. In \cite{B}, Beals showed:
\thm{4.1}{$S^0_{\Phi,\gvp}$ is a $\Psi$-subalgebra of  $\cL(L^2(\R^n))$.}

We now modify the hypoellipticity condition \ref{2.0} to

\dfn{4.2}{Let $a\in S^\mu_{\Phi,\gvp}$ for some $\mu\ge0$,
possibly matrix-valued. Assume
there exist constants $c,C>0$ such that for all $x,\xi\in\R^n$,
with $(\gvp\Phi)(x,\xi)\ge C$,
the spectrum of $a(x,\xi)$ lies outside  $\gL\cup\{|\gl|\le c\}$,
and, for $\gl\in\gL$,
\begin{eqnarray}\label{e4.2}
|\partial^\ga_\xi \partial^\gb_x  a(x,\xi)|\ |(a(x,\xi)-\gl)^{-1}|\le
c_{\ga\gb} \Phi^{-|\alpha|}\gvp^{-|\beta|}.
\end{eqnarray}
We define $\gO_{x,\xi}$ as before and note that \eqref{e4.2} extends to
$\C\setminus\gO_{x,\xi}$ for $(\gvp\Phi)(x,\xi)\ge C$.
}

While the $b_j$ are as before, we make a small change in the construction of $b^N$
in order to account for the new hypoellipticity condition:
Instead of the function $\gvp$
employed there, we use a zero-excision function $\tilde \gvp$ defined as follows:
We choose a function $\psi\in C^\infty(\R)$  with $\psi(t)=0$ for $t\le 1$ and
$\psi(t)=1$ for $t\ge 2$ and then let $\tilde\gvp(x,\xi)= \psi((\Phi\gvp)(x,\xi)/C)$
with the above constant $C$.

The results of \ref{2.10} - \ref{2.17} then follow in an analogous way,
replacing in Lemma \ref{2.10} and Remark \ref{2.12}
the space $S^0_{\rho,\gd}$ by $S^0_{\Phi,\gvp}$ and  $S^{m-N(\rho-\gd)}_{\rho,\gd}$
by $S^{\mu-(N,N)}_{\Phi,\gvp}$ with the above definition of $(K,k)$. Note that
the condition in Corollary \ref{2.14} can be fulfilled as a consequence of the
assumptions on $\mu$.

We then obtain the statements of Theorem \ref{3.01} (with $f(a)$ now in $S^0_{\Phi,\gvp}$),
Corollary \ref{3.02} and Theorem \ref{3.1} in the same way as before.
 
\section{Symbols in Weyl-Hörmander Classes}

Let $\gs$ be the canonical symplectic form on $\R^{2n}$ and  
$g$ a
$\gs$-temperate metric, see Hörmander \cite{H} for
details.
We denote by $|\cdot|_k^g$, $k\in \N_0$, the associated
seminorm system.
Let $m$ be a $(\gs,g)$-temperate weight function and $h,
H$ be defined by
$$h(x,\xi)^2=
\sup_{(y,\eta)}\frac{g_{(x,\xi)}(y,\eta)}{g^\gs_{(x,\xi)}(y,\eta)},\quad
H(x,\xi)^2=\sup_{(y,\eta)}\frac{g^\gs_{(x,\xi)}(y,\eta)}{g_{(x,\xi)}(y,\eta)}
.$$
where $g^\gs$ is the dual metric to $g$. It is required
that $h\le 1$.
Note that $h$ is a $(\gs,g)$-temperate  weight function and, by suitable 
regularization \cite[Section 2]{H}, $h$ can be assumed to be smooth.  
The symbol space $S(m,g)$ then consists of all functions
$p$ on $\R^n\times\R^n$
for which $|p|^g_k(x,\xi)/m(x,\xi)$ is bounded for each
$k$.
A symbol $p\in S(m,g)$ defines the Weyl pseudodifferential
operator $\op^wp:\cS\to\cS$ by
$$(\op^wp) u(x) =(2\pi)^{-n}\ \iint e^{i(x-y)\xi}
p((x+y)/2,\xi)u(y)\,dy\, d\xi.$$
For $p\in S(m_1,g)$ and $q\in S(m_2,g)$ we have $\op^w
p\circ\op^w q=\op^w r$ for
some $r\in S(m_1m_2,g)$, and the symbol seminorms for $r$
can be estimated in
terms of those for $p$ and $q$. We write $r=p\#^wq$.
Furthermore, the remainder
\begin{eqnarray*}
R_N(p,q)=
(p\#^wq)(x,\xi)-\sum_{j<N}
\frac{\sigma(\partial_x,\partial_\xi;\partial_y,\partial_\eta)^j}{(2i)^j j!}
p(x,\xi)q(y,\eta)|_{(y,\eta)=(x,\xi)}
\end{eqnarray*}
is an element of  $S(m_1m_2h^N,g)$.
The corresponding mapping is  continuous.

For $p$ in $S(1,g)$, $\op^wp$ is bounded on $L^2(\R^n)$,
so that
$S(1,g)$ can be considered a Fréchet sublgebra of
$\cL(L^2(\R^n))$ with the product $\#^w$.
It is not clear whether this algebra is
always spectrally invariant. This is, however, true under
mild restrictions, as shown by Bony \cite[Corollaire
4.4]{Bony}:

\thm{6.1}{$S(1,g)$ is a $\Psi$-subalgebra of
$\cL(L^2(\R^n))$ whenever $g$ is geodesically temperate
and of reinforced slowness $($`lenteur renforcée'$)$. }

Here, the metric $g$ is called geodesically temperate if
there exist constants $c>0$ and
$N\in \N$ such that for all $(x,\xi)$ and $(y,\eta)$
$$
\left(g^\gs_{(x,\xi)}(\cdot)\big/g^\gs_{(y,\eta)}(\cdot)\right)^{\pm1}
\le c\Big(1+d_\gs\big((x,\xi),(y,\eta)\big)\Big)^N,$$
where $d_{\gs}((x,\xi),(y,\eta))$ is the geodesic distance
of $(x,\xi)$ and
$(y,\eta)$ with respect to the metric $g^\gs$. Moreover,
$g$ is said to be of reinforced slowness, if there exists
a
constant $c>0$ such that
$$\left(g_{(x,\xi)}(\cdot)\big/g_{(y,\eta)}(\cdot)\right)^{\pm1}\le
c
~\text{ whenever }~g^\gs_{(x,\xi)}((x,\xi)-(y,\eta))\le
c^{-1}H^2(x,\xi).$$

\extra{6.2}{Hypoellipticity assumption}{Let $a$ be a complex-valued
symbol in $S(m,g)$, with $g$
geodesically temperate and of reinforced slowness and,
moreover,
\begin{eqnarray}\label{6.1.0}
c_m h^L\le m\le c_m^{-1}h^{-L}\ \text{for suitable } c_m,L>0.
\end{eqnarray}
Assume that there
exist constants $c,c'$ such that 
$a(x,\xi)$ lies outside $\gL\cup \{|\gl|\le c\}$ whenever
$h(x,\xi)\le c'$
and that, for $\gl\in\gL$ and suitable constants $c_k$,
\begin{eqnarray}\label{6.2.1}
|a|_k^g(x,\xi) \ |(a(x,\xi)-\gl)^{-1}|\le c_k,\quad
\gl\in\gL, k\in \N_0.
\end{eqnarray}
}

\extra{6.3}{Parametrices and inverses in the Weyl
calculus}{For $(x,\xi)$ with $h(x,\xi)\le c'$ we define
$\gO_{x,\xi}$ as before and note that
estimate \eqref{6.2.1} extends to $\C\setminus
\gO_{x,\xi}$.

Next we construct coarse parameter-dependent right and
left parametrices $b^N$
and $\tilde b^N$ to $a-\gl$ with respect to
the Weyl symbol product.
For $\gl\notin\gO_{x,\xi}$ and
all $(x,\xi)$ with $h(x,\xi)\le c'$ we determine
$b_j=b_j(x,\xi;\gl)$ iteratively as follows:
\begin{align*}
b_0(x,\xi;\gl)&= (a(x,\xi)-\gl)^{-1},\\
b_{j}(x,\xi;\gl)&= -b_0(x,\xi;\gl)\sum_{k+l=j;~
k>0}\frac{\gs^k}{(2i)^kk!}
a(x,\xi)b_l(y,\eta)|_{(y,\eta)=(x,\xi)},
\end{align*}
$j=1,2,\dots$, where for better legibility we wrote $\gs$
instead of $\gs(\partial_x,\partial_\xi;\partial_y,\partial_\eta)$.

We choose a smooth function $\psi$ on $\R$ with
$\psi(t)=1$ for $t<1/2$ and $\psi(t)=0$
for $t\ge1$ and let $\gvp(x,\xi)=\psi(h(x,\xi)/c')$ with
the above constant $c'$.
We then define
$$b^N(x,\xi;\gl)=\gvp(x,\xi) \sum_{j<N} b_j(x,\xi;\gl).$$
Next we recall that, for two smooth functions $p$, $q$,
\begin{align}\label{BN.1}
\begin{split}
|(\gs(\partial_x&,\partial_\xi;\partial_y,\partial_\eta)^j
p(x,\xi)q(y,\eta)|_{(y,\eta)=(x,\xi)}|_k^g\\
&\le (2n)^j\sum_{l=0}^k\binom kl
|p|^g_{j+l} |q|^g_{j+k-l}h^j,
\end{split}
\end{align}
see e.g.\ Buzano and Nicola, \cite[Lemma 2.3]{BN}.
It follows from \eqref{6.2.1} and \ref{2.4} that
$$|\gvp b_0(\gl)|^g_k(x,\xi)\le c_k|a(x,\xi)-\gl|^{-1}$$
for suitable constants $c_k$.
This shows that $\skp\gl \gvp b_0(\gl)$ is bounded in $S(1,g)$.
Moreover, we infer from \eqref{6.2.1}, \eqref{BN.1},
and the iteration process that
\begin{eqnarray}\label{bj}
|\gvp b_j(\gl)|_k^g(x,\xi) \le c_k|a(x,\xi)-\lambda|h(x,\xi)^j
\end{eqnarray}
with (different) constants $c_k$, uniformly for all $x,\xi$
and $\lambda\not\in\Omega_{x,\xi}$.}

\lemma{6.4}{Let $r^N(\gl)= (a-\gl)\#^wb^N(\gl) -1$. Then
$\gl\mapsto \skp\gl r^N(\gl)$ is a bounded continuous map
from $\gL$ to
$S(mh^N,g)$.}

\begin{proof}
The proof is analogous to that of Lemma \ref{2.10}. We
write $r_N(\lambda)$ as
\begin{equation}\label{6.2.2}
 R_N\left(a-\lambda,b^N(\lambda)\right)    
+\Big(\sum_{j,k<N}\frac{\gs^k}{(2i)^kk!}
(a(x,\xi)-\lambda)(\gvp b_j)(y,\eta;\gl)|_{(y,\eta)=(x,\xi)}-1\Big).
\end{equation}
 Then $\skp\gl
R_N\left(a-\lambda,b^N(\lambda)\right)=R_N\left(a,\skp\gl
b^N(\lambda)\right)$
 is bounded and continuous in $\lambda$ with values in
$S(mh^N,g)$.
Next we consider the second summand. We split the summation into those terms,
where $j+k<N$ and those, where $j+k\ge N$.
By construction, the former sum equals $\gvp$.
The difference $\gvp-1$ belongs
to $S(h^M,g)$ for any $M$, and hence to $S(mh^N,g)$.
In the latter, we have $j,k\not=0$ and can therefore
replace $a-\gl$ by $a$. We conclude from \eqref{6.2.1}, \eqref{BN.1} and \eqref{bj},
that, even after multiplication by $\skp\gl$,  
the $S(mh^N,g)$-seminorms of all these terms are finite.
Continuity follows from the resolvent identity as before. 
\end{proof}

Similarly, we construct a left parametrix  $\tilde b^N$ with remainder $\tilde r^N$ such that
$\gl\mapsto \skp\gl \tilde r^N(\gl)$ is a bounded continuous map from $\gL$ to $S(mh^N,g)$.

Since, by assumption, $S(mh^N,g)\hookrightarrow S(1,g)$ for large $N$,
the results of \ref{2.14} - \ref{2.17} then follow as above.

\extra{6.10}{Functional calculus}{We
obtain the statements of Theorem \ref{3.01} with $f(a)$ now in $S(1,g)$,
of Corollary \ref{3.02} and Theorem \ref{3.1} in the same way as before.
}

\begin{small}
\bibliographystyle{amsalpha}

\end{small}

\end{document}